\documentclass[12pt]{article}
\usepackage{amsmath}
\usepackage{amssymb}
\newcommand{\Ann}{\mbox{Ann}}

\newcommand{\Ass}{\mbox{Ass}}
\newcommand{\Supp}{\mbox{Supp}}

\newcommand{\depth}{\mbox{depth}}

\newcommand{\Min}{\mbox{Min}}

\begin{document}

\begin{center}
\large{\bf Bounds for numbers of generators for a class of
submodules of a finitely generated module}
\end{center}

\vspace{.2in}

\begin{center}
{\bf Tirdad Sharif and Siamak Yassemi\footnote{This research was
supported in part by a grant from IPM.}}
\end{center}

\begin{center}
{\it Department of Mathematics, University of Tehran\\and\\ Institute for Studies in
  Theoretical Physics and Mathematics (IPM).}
\end{center}

\vspace{.3in}

\begin{abstract}
We determine a bound for numbers of minimal generators for a
class of submodules of a finitely generated module $M$ with
dimension $d\ge 1$.
\end{abstract}

\vspace{.2in}



\vspace{.3in}

\baselineskip=18pt

\noindent{\large\bf 0. Introduction}

\vspace{.2in}

Throughout this note the ring $R$ is commutative Noetherian with
non--zero identity and $M$ is a finite (that is, finitely
generated) $R$--module. We use the notations $\mu(K)$ and
$\ell(K)$ respectively for the minimal number of generators and
the length of an $R$--module $K$. It is a well-known fact that a
local ring is of dimension at most one if and only if there exists
a non--negative integer $n$ such that $\mu(I)\le n$ for all ideals
$I$ of $R$, cf. [{\bf S}; Theorem 1.2, chapter 3]. For a
zero--dimensional local ring, J. Watanabe found the following
theorem: for any ideal $I$ of $R$ and any element $x\in\frak m$,
we have $\mu(I)\le\ell(R/xR)$, cf. [{\bf W}]. To find uniform
bounds in rings of dimension greater than one, we must restrict
ourselves to appropriate subclasses of the class of all ideals in
the ring. In [{\bf G}], Gottlieb proved one such result: Let $R$
be a local ring of dim $d \ge 1$. Let $q$ be an ideal generated by
a system of parameters. Then $\mu(I)\le\ell(R/q)$ for all $I$ such
that depth $R/I \ge d-1$.

In [{\bf M}], Matsumura proved the following theorem for finite
modules:


\noindent{\bf Theorem A.}([{\bf M}; Theorem 2]) Let $(R,\frak m)$
be a local ring and let $M$ be a finite $R$--module of dimension
at most one. Let $N$ be a submodule of $M$ and let $x\in\frak m$.
Then $\mu(N)\le\ell(M/xM)$.

The aim of this paper is to obtain a uniform bound for a certain
class of submodules from the following theorem: Let $(R,\frak m)$
be a local ring, let $M$ be a finite $R$--module of dimension
$d\ge 1$ and let $\frak q$ be an ideal of $R$ generated by a
system of parameters on $M$. Let $N$ be a submodule of $M$ with
$\depth M/N\ge d-1$. Then $\ell(N/\frak qN)\le\ell(M/\frak qM)$.

\vspace{.3in}

\noindent {\large\bf 1. Main results}

\vspace{.2in}

First of all we bring a generalization of [{\bf S}; Theorem 1.2,
chapter 3] for modules.

\vspace{.2in}

\noindent{\bf Theorem 1.} Let $(R,\frak m)$ be a local ring and
let $M$ be a finite $R$--module of dimension $d$. Then $d\le 1$ if
and only if there is a non--negative integer $t$ such that
$\mu(N)\le t$ for every submodule $N$ of $M$.

\vspace{.1in}

\noindent{\it Proof.} Assume that for every submodule $N$ of $M$,
we have $\mu(N)\le t$. By the Hilbert Samuel Theorem, cf. [{\bf
S}; Theorem 2.1' page 5], we have for sufficiently large $n$,
$\ell(\frak m^nM/\frak m^{n+1}M)=\mu(\frak m^nM)$ is a polynomial
in $n$ of degree d-1. Hence $d\le 1$. The converse is obvious from
[{\bf M}; Theorem 2].

\vspace{.2in}

 The following theorem is a generalization of Theorem B because
$\mu(N)\le\ell(N/xN)$ for any $x\in\frak m$.

\vspace{.1in}

\noindent{\bf Theorem 2.} Let $(R,\frak m)$ be a local ring and
let $M$ be a finite $R$--module with $\dim M\le 1$. Then for any
$x\in\frak m$ and $N$ a submodule of $M$, we have
$\ell(N/xN)\le\ell(M/xM)$.

\vspace{.1in}

\noindent{\it Proof.} If $\ell(M/xM)=\infty$, the inequality is
clear. Therefore assume that $\ell(M/xM)<\infty$. Since $\dim
M/xM=0$, by using the fact that $\sqrt{\frak q R + \Ann(M)} =
\sqrt{(\frak q M:_RM)}$ for any ideal $\frak q$ of $R$, we have
that $\dim R/(xR+\Ann(M))=0$ ``and hence there'' exists
$s\in\mathbb N$ such that $\frak m^s\subseteq xR+ \Ann(M)$. Thus
$\frak m^sN\subseteq xN$. Therefore $\ell(N/xN)$ is finite. Now we
consider two cases:

Case 1. $\ell(M/N)<\infty$. Since $xM/xN\cong M/(xN:_Mx)$ is
isomorphic to a quotient of $M/N$, we have that
$\ell(M/N)=\ell((xN:_Mx)/N)+\ell(xM/xN)$. Consider the exact
sequences
$$0\to N/xN\to M/xN\to M/N\to 0,$$
$$0\to xM/xN\to M/xN\to M/xM\to 0.$$

We have
$\ell(M/xM)-\ell(N/xN)=\ell(M/N)-\ell(xM/xN)=\ell(xN:_Mx/N)\ge 0$.
Now the inequality holds.

Case 2. $\ell(M/N)=\infty$. We can reduce this case to Case 1. By
the Artin--Rees theorem there exists $t\in\mathbb N$ such that
$\frak m^tM\cap N\subseteq\frak m^sN\subseteq xN$. Therefore
$$N/xN=N/(xN+(\frak m^tM\cap N))=N/(N\cap(\frak m^tM+xN))=(N+\frak
m^tM)/(\frak m^tM+xN).$$ Set $K=N+\frak m^tM$. We claim that
$xN+\frak m^tM\supsetneq xK$. If not, we have $\frak
m^t(M/xN)=x\frak m^t(M/xN)$ and so by Nakayama's Lemma $\frak
m^tM\subseteq N$. Thus $\dim M/N=0$, a contradiction. Therefore
$\ell(N/xN)=\ell(K/(xN+\frak m^tM)<\ell(K/xK)$. Now the assertion
follows from Case 1 and the fact that $\dim M/K=0$.

\vspace{.2in}

\noindent{\bf Theorem 3.} Let $(R,\frak m)$ be a local ring and
let $M$ be a finite $R$--module of dimension $d\ge 1$. Let $N$ be
a submodule of $M$ with $\depth M/N\ge d-1$. Then for any ideal
$\frak q$ that is generated by a system of parameters on $M$,
$\ell(N/\frak qN)\le\ell(M/\frak qM)$.

\vspace{.2in}

\noindent{\it Proof.} We prove it by induction on $d$. For $d=1$
the assertion follows from Theorem 2. Suppose that $d>1$,
$\Min(\Ann M)=\{\frak p_1,\frak p_2,\dots,\frak p_t\}$ and $\Ass
(M/N)=\{\frak p_{t+1},\dots,\frak p_n\}$. Since $\dim M>0$ and
$\depth M/N>0$ we have that $\frak p_i\neq \frak m$ for each $
i\le n$. Let $\frak q$ be an ideal of $R$ that is generated by a
system of parameters on $M$. Since $\Supp M/\frak qM=\{\frak m\}$
we have that $\frak q\nsubseteq\frak q\frak
m\cup(\cup_{i=1}^n\frak p_i)$. Choose $x\in\frak q\backslash\frak
q\frak m\cup(\cup_{i=1}^n\frak p_i)$. Set $\overline{T}=T/xT$
(where $T$ is one of $R$ or $M$) and $\widetilde{R}=R/\Ann M$.
Then we have
$$\begin{array}{rl}
\dim_{\overline{R}}\overline{M} &=\,
\dim\overline{R}/\Ann_{\overline{R}}\overline{M}=\dim
R/\Ann_R\overline{M}\\
&=\, \dim R/(\Ann_RM+xR)=\dim\widetilde{R}/x\widetilde{R}\\ & =\,
\dim\widetilde{R}-1=\dim M-1.
\end{array}$$

On the other hand we have that
$\depth_{\overline{R}}\overline{M}/\overline{\overline{N}}=
\depth_R\overline{M}/\overline{\overline{N}}$ where
$\overline{\overline{N}}=(N+xM)/xM$, cf. [{\bf BH}; Exercise
1.2.26] and hence
$$\begin{array}{rl}
\depth_{\overline{R}}\overline{M}/\overline{\overline{N}} & =\,
\depth_RM/(N+xM)=\depth (M/N)/x(M/N)\\ &=\, \depth_RM/N-1\ge d-2.
\end{array}$$

Since $x\notin\frak q\frak m$ we can extend $x$ to a minimal
generating set for $\frak q$ and hence $\frak Q=\frak q/xR$ is
generated by a system of parameters on $\overline{M}$. Thus by the
induction hypothesis
$\ell_{\overline{R}}(\overline{\overline{N}}/\overline{\frak
q}\overline{\overline{N}})\le\ell_{\overline{R}}(\overline{M}/\frak
Q\overline{M})$. It is easy to see that $\frak Q\overline{M}=\frak
qM/xM$. Therefore $\ell_{\overline{R}}(\overline{M}/\frak
Q\overline{M})=\ell_R(M/\frak qM)$ and $\frak
Q\overline{\overline{N}}=(\frak qN+xM)/xM$. Also,
$\ell_{\overline{R}}(\overline{\overline{N}}/\frak Q
\overline{\overline{N}})=\ell_R((N+xM)/\frak (qN+xM))$. Set
$\ell_R(N/\frak qN)=t$. Then there is a chain $\frak q N = N_0
\subset \dots\subset N_t = N$ with $N_i/N_{i-1} \cong R/\frak m$
for $1 \le i \le t$. We claim that the chain $N_0+xM \subseteq
\dots\subseteq N_t+xM$ is also strictly increasing (necessarily
with simple factors). If not, we have $N_i \subseteq N_{i-1} +xM$
for some index $i\ge 1$. Let $a\in N_i$. Let $a\in N_i$. Then
$a\in N_{i-1}+xM$ and hence $a=b+xm$ for some $b\in N_{i-1}$ and
$m\in M$. Since $xm\in N$ and $x$ is not a Zero--divisor on $M/N$
we have that $m\in N$. Therefore $N_i\subseteq N_{i-1}+xN\subseteq
N_{i-1}+\frak qN\subseteq N_{i-1}$. This is a contradiction. Thus
$\ell_{\overline{R}}(\overline{N}/\overline{\frak
q}\overline{N})=\ell(N/\frak qN)$ and the assertion holds.

\vspace{.2in}

\noindent{\bf Corollary 4 (Main Theorem).} Let $(R,\frak m)$ be a
local ring and $M$ a finite $R$-module of dimension $d \ge 1$. Let
$\frak q$ be an ideal generated by a system of parameters on $M$,
and put $c = \ell(M/\frak q M)$. Then, for any submodule $N$ of
$M$ with $\depth(M/N) \ge d-1$, we have $\mu(N) \le c$.

\vspace{.2in}

\noindent{\it Proof.} Since $\frak q\subseteq\frak m$ we have
that $\mu(N)=\ell(N/\frak mN)\le\ell(N/\frak qN)\le\ell(M/\frak
qM)$.

\vspace{.2in}

\noindent{\bf Corollary 5.} Let $(R,\frak m)$ be a local ring and
let $M$ be a finite Cohen--Macaulay $R$--module of dimension $d\ge
1$. Let $\frak q$ be an ideal generated by a system of parameters
on $M$, and put $c=\ell(M/\frak qM)$. Let $N$ be a Cohen--Macaulay
submodule of $M$ with $\dim N=d$. Then $\mu(N)\le c$, and every
$M$--sequence is an $N$--sequence.

\vspace{.2in}

\noindent{\it Proof.} It follows from [{\bf BH}; Proposition
1.2.9] that $\depth(M/N) \ge d-1$. The inequality now follows from
Corollary 4. Assume that $x=x_1, x_2, \dots, x_r$ is an
$M$--sequence. Since $x$ is a part of a system of parameters
$\frak q$ on $M$, cf. [{\bf BH}; Theorem 2.1.2], we have that
$\ell(N/\frak qN)\le\ell(M/\frak qM)$ by Theorem 3. Therefore
$\frak q$ is a system of parameters on $N$. Thus $\frak q$ is
generated by an $N$--sequence, cf. [{\bf BH}; Theorem 2.1.2], and
hence $x$ is an $N$--sequence.

\vspace{.2in}

\noindent{\bf Corollary 6.} Let $R$ be a Cohen-Macaulay ring, let
$\frak q$ be a parameter ideal of $R$, and let $t = \ell(R/\frak
q)$. Let $M$ be a maximal Cohen-Macaulay module of rank $r$ (that
means $M_P$ is free of constant rank $r$ for $P\in \Ass(R)$). Then
$\mu(M)\le rt$.

\vspace{.1in}

\noindent{\it Proof.} It is easy to see that the module $M$ is
isomorphic to a submodule of $R^r$. Now apply Corollary 5.

\vspace{.3in}




\baselineskip=16pt

\begin{center}
\large{\bf References}
\end{center}

\vspace{.2in}

\begin{verse}

[BH] W. Bruns and J. Herzog, {\em Cohen--Macaulay rings}, Cambridge
University Press, Cambridge 1993.



[G] C. Gottlieb, {\em Bounding the number of generators for a
class of ideals in local rings}, Comm. in Algebra, {\bf 23}
(1995), 1499--1502.


[M] H. Matsumura, {\em Some new results on numbers of generators
of ideals in local rings}, Topics in algebra, Part 2 (Warsaw,
1988), 157--161, Banach Center Publ., 26, Part 2, PWN, Warsaw,
1990.

[S] J. D. Sally, {\em Numbers of generators of ideals in local
rings}. Marcel Dekker, Inc., New York-Basel, 1978.

[W] J. Watanabe, {\em m--full ideals}, Nagoya Math. J. {\bf 106}
(1987), 101--111.

\end{verse}

\end{document}